# Model Averaging from a Geometrical Perspective

C. L. Winter[1] and Doug Nychka[2]

May 21, 2008

Abstract. Given a collection of computational models that all estimate values of the same natural process, we compare the performance of the average of the collection to the individual member whose estimates are nearest a given set of observations. Performance is the ability of a model, or average, to reproduce a sequence of observations of the process. We identify a condition that determines if a single model performs better than the average. That result also yields a necessary condition for when the average performs better than any individual model. We also give sharp bounds for the performance of the average on a given interval. Since the observation interval is fixed, performance is evaluated in a vector space, and we can add intuition to our results by explaining them geometrically. We conclude with some comments on directions statistical tests of performance might take.

**Introduction**

Computational models are often used to estimate states, $Y_t$, of a given physical process. In many cases, several alternative models are available to estimate $Y_t$. Good examples are the large number of climate models used by the IPCC to estimate future impacts of increased carbon in the atmosphere (ref's) and hydrogeology (ref's). In general, there are $m = 1, 2, \ldots, M$ models, each producing estimates $X_{m,t}$ of $Y_t$. The question naturally arises how to estimate $Y_t$ in light of the alternatives. One common approach is to approximate $Y_t$ with the weighted average of the $X_{m,t}$,

$$\bar{X}_t = \sum_m w_m X_{m,t}, \quad \sum_m w_m = 1, \quad w_m \geq 0 \text{ for all } m. \tag{1}$$

It has been observed that $\bar{X}_t$ sometimes produces "better" estimates of $Y_t$ than any individual model [1]-[9].

After defining "better" and "performance" in terms of the behavior of $Y_t$, $X_{m,t}$, and $\bar{X}_t$ on finite time intervals $t \in [t_0, t_T]$, we identify conditions that determine when $\bar{X}_t$ performs better (or not) over a given time interval than any individual model over that same interval. Our basic result is a necessary condition on the collinearity of models that must be met for the average to perform better. We also give a condition that implies the best individual model will outperform the average, and we put sharp bounds on the performance of the average. Although our basic outlook in this paper is statistical, our methods and results are geometrical and specific to a given interval; our results are "pre-statistical" in a sense that will become clearer as we go along. We make some general comments on model selection in later sections of this paper, and we apply our results to global climate models. We conclude by indicating a path for converting our results to statistical tests.

---

[1] National Center for Atmospheric Research, P.O. Box 3000, Boulder, CO 80307-3000, lwinter@ucar.edu.
[2] National Center for Atmospheric Research, P.O. Box 3000, Boulder, CO 80307-3000, nychka@ucar.edu.

**Model Performance on $[t_0, t_T]$**

We define a model's performance in terms of a measure over the set of all finite time intervals, $[t_0, t_T]$. In recognition of the discrete nature of time in computer simulations, we write $[t_0, t_T] \equiv \{t_0, t_0 + 1, \ldots, t_0 + T\}$ and we use corresponding discrete operators, but the discussion can be translated into continuous domains if necessary.

The performance of a model over a given interval $[t_0, t_T]$ can be evaluated by calculating its average (over time) mean-squared departure from known data, $Y_t$, $t \in [t_0, t_T]$,

$$S_m^2(t_0, t_T) = \frac{1}{T} \sum_{t=t_0}^{t_T} (X_{m,t} - Y_t)^2 . \tag{2}$$

Model $m$ outperforms model $m'$ on $[t_0, t_T]$ if $S_m^2(t_0, t_1) < S_{m'}^2(t_0, t_1)$. Other performance measures can be defined, but (2) is common and serves as a benchmark.

Similar to (2), the average model's performance on $[t_0, t_T]$ is

$$S^2(t_0, t_T) = \frac{1}{T} \sum_t (\overline{X}_t - Y_t)^2$$

$$= \frac{1}{T} \sum_t (\sum_m w_m X_{m,t} - Y_t)^2$$

$$= \sum_m w_m^2 S_m^2(t_0, t_T) + \sum_{m \neq m'} w_m w_{m'} R_{m.m'}(t_0, t_T). \tag{3}$$

The *correspondence* between models,

$$R_{m,m'}(t_0, t_T) = \frac{1}{T} \sum_t (X_{m,t} - Y_t)(X_{m',t} - Y_t), \tag{4}$$

measures their mutual agreement over $[t_0, t_T]$.

The statistical performance of models can be evaluated by calculating quantities like

$$\text{Pr}ob[S_m(t', t'') < S_{m'}(t', t'') \text{ for any } [t', t''] \mid \{S_m(t_0, t_T), S_{m'}(t_0, t_T)\} \text{ for some } \{[t_0, t_T]\}].$$

Ultimately we are interested in

(*) $\operatorname{Prob}[S(t',t'') < S_{m'}(t',t'')$ for all $m$ and any $[t',t''] | S_m(t_0,t_T), S_{m'}(t_0,t_T)$ for $[t_0,t_T]]$,

but our results in this paper, are limited to performance over a given interval $[t_0,t_T]$.

Since we limit ourselves to given $[t_0,t_T]$, we simplify our notation by dropping dependence on $[t_0,t_T]$ from expressions like $S_m$ and $R_{m,m'}$, and we just write $t = 1,2,\ldots,T$ for $t \in [t_0,t_T]$. Now we can easily interpret (1)-(3) in terms of a $T$-dimensional vector space, $\Re^T$. First, let $\vec{Z}_m = (X_{m,1} - Y_1, \ldots, X_{m,T} - Y_T)$ and $\vec{Z} = (\bar{Z}_1, \ldots, \bar{Z}_T)$. The squared departures are obviously equivalent to the $L^2$ distances

$$S_m^2 = \frac{1}{T}\|\vec{Z}_m\|^2, \tag{2*}$$

$$S^2 = \frac{1}{T}\|\vec{Z}\|^2, \text{ and} \tag{3*}$$

$$R_{m,m'} = S_m S_{m'} \cos\theta_{m,m'} \tag{4*}$$

with $\theta_{m,m'}$ the angle between $\vec{Z}_m$ and $\vec{Z}_{m'}$.

**Results**

Our basic result (Result 3 below) is that individual models cannot all correspond strongly on $[t_0,t_T]$ when the average performs better than any individual. We obtain it as the corollary to a preliminary result about the opposite case when at least one individual model performs better than the average. The best performer among a collection of models is $\vec{Z}_{\min}$ corresponding to $S_{\min}^2 = \min_m\{S_m^2\}$.

It is easy to find a criterion for when the best model outperforms the average on given $[t_0,t_T]$, i.e., when $S^2 > S_{\min}^2$.

Result 1. In the context of (1)-(4) and the discussion surrounding them, if $R_{m,m'} > S_{\min}^2$ for all $m$, $m'$, then $S^2 > S_{\min}^2$.

Proof. Clearly, $S^2 - S_{\min}^2 = \sum_m w_m^2(S_m^2 - S_{\min}^2) + \sum\sum_{m \neq m'} w_m w_{m'}(R_{m,m'} - S_{\min}^2)$. As $\sum_m w_m^2(S_m^2 - S_{\min}^2) \geq 0$, we only need $\sum\sum_{m \neq m'} w_m w_{m'}(R_{m,m'} - S_{\min}^2) \geq 0$. That is obvious because i) $w_m \geq 0$ for all $m$ and ii) by hypothesis $R_{m,m'} - S_{\min}^2 \geq 0$ for all $m$, $m'$. □

Since the hypothesis $R_{m,m'} \geq S^2_{min}$, is equivalent to $\cos\theta_{m,m'} > \frac{S^2_{min}}{S_m S_{m'}}$ for all $m$, $m'$, we can restate the Result 1 in terms of the geometry of $\mathfrak{R}^T$.

Result 2. In the context of (1)-(4), if $\cos\theta_{m,m'} > \frac{S^2_{min}}{S_m S_{m'}}$ for all $m$, $m'$, then $S^2 > S^2_{min}$.

Now we are ready for our main result. It follows easily from Results 1-2 that models cannot all correspond highly correlated on $[t_0, t_T]$ if the average is to perform better than the best model.

Result 3. In the context of (1)-(4), if $S^2 < S^2_{min}$, then $\cos\theta_{m,m'} < \frac{S^2_{min}}{S_m S_{m'}}$ for some $m$, $m'$.

Remark. It is interesting to look at special cases in the light of these facts.

  i) If all models are about equally good on $[t_0, t_T]$, i.e., $S^2_{min} \sim S^2_m$, but do not correspond highly, i.e., $\cos\theta_{m,m'} \ll 1$, then Results 1-2 will not usually apply because usually $\frac{S^2_{min}}{S_m S_{m'}} \sim 1 \gg \cos\theta_{m,m'}$.
  ii) If one model is much better than the rest, then $S^2_{min} \ll S^2_m$ for all $m \neq min$. Results 1-2 will usually apply in cases of positive correspondence between models, since then $\cos\theta_{m,m'} > \frac{S^2_{min}}{S_m S_{m'}}$.

Result 4. Using the Schwartz inequality,

$$0 \leq S^2 \leq \sum_m w_m^2 S_m^2 + \sum_{m \neq m'} w_m w_{m'} S_m S_{m'} = (\sum_m w_m S_m)^2. \tag{5}$$

These bounds are sharp. The upper bound is reached when all models are the same or for that matter, are just in perfect correspondence, i.e., $\cos\theta_{m,m'} = 1$ for all $m$, $m'$.

**Geometry of Model Evaluation**
The usual case in model evaluation is that a limited number of observations, $\vec{Y}$, are available for $[t_0, t_T]$, a restricted interval of time [1]. Often $[t_0, t_T] = [t_0, t_N] \cup [t_{N+1}, t_T]$ is split into two subintervals, a calibration interval $C = [t_0, t_N]$ used to parameterize models and a validation interval $V = [t_{N+1}, t_T]$ used to evaluate the model's performance. Models are evaluated by comparing $S^2(V)$ and $\{S_m^2(V)\}$, or comparable metrics. The comparisons are not statistical, but rather are restricted to the given $V$, so Results 1-3 apply. In many cases (possibly, most), different models will yield $S^2_{min}$ for different intervals; hence the need for statistical formulations based on (*).

The geometry of the space $\Re^T$ defined on a given validation interval $V$ clarifies our results [Figures 1-2]. In the Figures, $T = 3$ for convenience of illustration, so each $\vec{Z}_m = (Z_{m,1}, Z_{m,2}, Z_{m,3})$. The thin vectors represent the performance of different models, while the thick vector is the average, $\bar{Z}$. Figure 1 shows a typical case when $S^2 \geq S^2_{min}$ because model performance is roughly collinear. The models all vary similarly about $\vec{Y}$, as evidenced by their collinearity, but one, $\vec{Z}_{min}$, performs much better than the others. The result of averaging is to "stretch" $\vec{Z}_{min}$ in the general direction of the other models. Figure 2 is a typical case when $S^2 \leq S^2_{min}$. The role of the anti-correspondence requirement -- $\cos\theta_{m,m'} < \dfrac{S^2_{min}}{S_m S_{m'}}$ for some $m$, $m'$ -- is clear: $S^2_{min}$ exceeds $S^2$ because the two sets of vectors "pull" against each other to produce $\bar{Z}$.

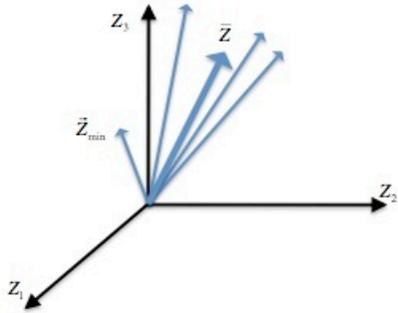
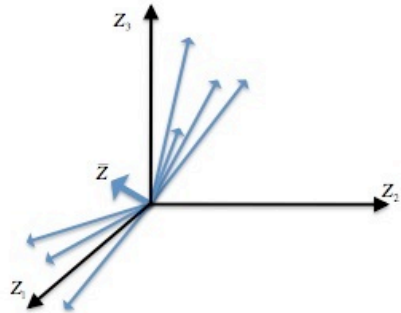

Figure 1 – Example of Result 1. Collinearity places $\bar{Z}$ beyond $\vec{Z}_{min}$.

Figure 2 – Example of Result 3. Weak correspondence allows $S^2 \leq S^2_{min}$.

**Global Climate Model Averages**
The IPCC has considered cases where the weights in (1) are uniform: $w_m = 1/m$ for all $m$.

**Relation to Model Evaluation**
The obvious way to make $S^2$ small is to have good models, $\vec{X}_m \approx \vec{Y}$, for all $m$ and all $[t_0, t_T]$. The usual case, however, is that we have only one, or at most a few, evaluation intervals, so currently evaluations are made on an ad hoc basis that is specific to the given interval(s). No statistical tests are available.

When some models perform well, but others don't, the set of models can be pre-screened according to a criterion like $\frac{S_m^2}{\|\vec{Y}\|^2} \ll 1$ [1,6]. Weighting good models more heavily than poor ones is a related way to reduce (3), but it begs the question, why include poor models in (1)? There may be a good answer, but it will be based on reasons that go beyond performance.

The correspondence between models, $R_{m,m'}$, can complicate model selection when all models are about equally bad, i.e., when

$$S_m^2 \approx S_{m'}^2 \text{ and } \frac{S_m^2}{\|\vec{Y}\|^2} \gg 1 \text{ for all } m, m'. \tag{6}$$

In this case, (3) suggests a quite different approach to selecting models: choose them so they're as anti-correlated as possible, i.e., choose $\sum_{m \neq m'} w_m w_{m'} R_{m,m'}$ as near $-\frac{1}{2}\sum_m w_m^2 S_m^2$ as possible. The upshot will be something like Figure 2.

**Discussion**
Evaluation of model averages should be based on statistics. The geometric view taken in this paper suggests an approach.